\documentclass[12pt]{article}
\usepackage{epsfig}
\usepackage{graphicx,subfigure,amsmath,url}
\include{lmacros}

\newcommand{\Ignore}[1]{}

\newcommand{\A}{\mbox{${\bf A}$}}
\newcommand{\B}{\mbox{${\bf B}$}}

\newcommand{\am}{\mbox{${\bf a}$}}
\newcommand{\bm}{\mbox{${\bf b}$}}

\newcommand{\ea}{\tilde{\bf a}}
\newcommand{\eb}{\tilde{\bf b}}
\newcommand{\ec}{\tilde{\bf c}}

\newcommand{\setRn}{R^n}
\newcommand{\covellipse}[2]{{\bf {\cal E}}(#1,#2)}
\newcommand{\ellipseTriple}{ \{A_i, b_i, c_i\} }

\title{\Large\bf On the Equivalence of the\\ General Covariance Union (GCU) and\\
       Minimum Enclosing Ellipsoid (MEE) Problems}

\author{{\bf Ottmar Bochardt} and {\bf Jeffrey Uhlmann}\\
Department of Computer Science, 201 EBW,\\
University of Missouri, Columbia, MO 65211\\
Phone: 573 884-2129, Fax: 573 882-8318}
\date{}
\begin{document}
\maketitle 

\begin{abstract}
In this paper we describe General Covariance Union (GCU) and 
show that solutions to GCU and the 
Minimum Enclosing Ellipsoid (MEE) problems are
equivalent. This is a surprising result because GCU is defined
over positive semidefinite (PSD) matrices with statistical 
interpretations while MEE involves PSD matrices with geometric 
interpretations. Their equivalence establishes an intersection 
between the seemingly disparate methodologies of covariance-based 
(e.g., Kalman) filtering and bounded region approaches to data fusion. 
\end{abstract}

\section{Introduction}

Positive semidefinite matrices are often chosen to represent
uncertainty in state estimation and control algorithms because
of their special linear-algebraic properties. The Kalman filter,
for example, uses PSD matrices to reprepresent covariance upper 
bounds on the second central moments of unknown probability
distributions relating to the state of a system of interest. 
Bounded region filters, by contrast, use PSD matrices to define 
ellipsoidal regions which bound the state of the system of 
interest. 

The choice between representing uncertainty with covariance
upper bounds versus ellipsoidal bounded regions leads to
very different data fusion algorithms with very different 
filtering and control properties. The statistical interpretation
associated with covariance matrices leads to use of the
Kalman fusion equations in the case of independent estimates or
Covariance Intersection (CI) when the statistical relationships 
among estimates to be fused cannot be established\cite{phdthesis:uhlmann95}. 
Under their
respective assumptions, Kalman and CI yield a fused estimate
with the smallest possible covariance that is guaranteed to be
an upper bound on the second central moment of the unknown
probability distribution defining the state of the system.
A bounded region filter, on the other hand, determines the 
minimum-size ellipsoid that bounds the intersection of the
given ellipsoids which are assumed to bound the true state
of the system. In summary, one framework is statistical while
the other is geometric.

Despite the use of PSD matrices and the exploitation of 
linearity properties, covariance and ellipsoidal frameworks 
employ very different mathematical techniques and have almost 
completely disjoint literatures. (This relationship also holds true
more generally between statistics and computational geometry.)
A consequence of this fact is that the sophisticated tools 
developed separately within each framework seem to have no 
applicability within the other framework. In this paper we
make progress toward bridging this gap by establishing a 
surprising equivalence result. 

The structure of this paper is as follows: The next section
provides background for the Kalman filter, Covariance 
Intersection (CI), and Covariance Union (CU). We then 
present General Covariance Union (GCU) and show that its
solution is equivalent to determining the Minimum Enclosing
Ellipsoid (MEE) of the ellipsoids defined by its covariance
matrix arguments. In other words, GCU and MEE are equivalent.
We end with a discussion of the implications of this 
unexpected result.

\section{The Covariance Upper Bound Framework for Data Fusion}

An estimate of the squared error, 
or {\em covariance}, associated 
with measurements from a particular sensor
{\em can} be modeled empirically by examining measurements 
taken of reference objects whose true states are known.  This
permits an error covariance matrix to be associated with 
subsequent measurements of objects whose true states are
not known. For example, the measured position of an object 
in two dimensions can be represented as a vector ${\bf a}$ 
consisting of the object's estimated mean position, e.g., 
${\bf a} = [{\bf x},{\bf y}]^T$, and an error covariance matrix 
${\bf A}$ that expresses the uncertainty associated with the 
estimated mean.  If the error in the estimated mean vector is 
denoted as $\ea$, then the error covariance
matrix is an estimate of the expected squared error, E$[\ea\ea^T]$. 

Ideally, ${\bf A}$ would equal E$[\ea\ea^T]$, but this is not
generally possible to achieve in practice because measurement and
process models are never perfect. To accommodate the effect of
model errors, prediction and measurement covariances are typically
overestimated so as to avoid underestimating the actual squared errors.
In other words, a more conservative overestimate of errors is 
deemed preferable to underestimating the errors. One of many reasons
for this preference is that it avoids the consequences of spuriously
small errors causing a covariance matrix to become singular
or numerically unstable. 

Formally, a mean and covariance estimate is said to be
consistent (or conservative) if and only if 
${\bf A}$ - E$[\ea\ea^T]$ is positive semidefinite 
(i.e., has no negative eigenvalues):
\begin{equation}
{\bf A} ~ \succeq ~ \mbox{E}[\ea\ea^T] 
\end{equation}

The statistical properties associated with a mean vector and its
associated covariance upper bound are completely defined by the
definion of consistency: all that can be said about a 
mean and covariance pair, $({\bf a}, {\bf A})$, is that 
the covariance matrix, ${\bf A})$, is greater than the expected
squared error in its associated mean, ${\bf a})$. 

Having established a rigorous definition of what constitutes
a {\em consistent} or {\em conservative} estimate, it is possible
to certify the performance of the Kalman filter and 
Covariance Intersection. 

\subsection{The Kalman Filter}

Given two mean and covariance estimates $({\bf a}, {\bf A})$
and $({\bf b}, {\bf B})$, the data fusion problem of interest
in this paper consists of
determining a fused estimate $({\bf c}, {\bf C})$ that is
guaranteed to be consistent and summarizes the information
in the two estimates with error (in terms of the size of
${\bf C}$) that is less than or equal
to that of either estimate.  If the two estimates are 
consistent and presumed to be statistically independent,
then a joint estimate can be constructed as:
\begin{equation}
\left( \left[\begin{array}{c} \am \\ \bm \end{array} \right],
\left[\begin{array}{cc} 
   \A  & {\bf 0} \\
   {\bf 0} & \B   
\end{array} \right] \right).
\end{equation}
Letting $\ea$ and $\eb$ denote the errors in the
respective mean estimates, the key property of the 
joint covariance estimate is that
it satisfies:
\begin{equation}
\left[\begin{array}{cc}
   {\bf A} & {\bf 0}\\
   {\bf 0} & {\bf B}
\end{array}\right]        ~ \succeq ~
\left[\begin{array}{cc}
   \mbox{E}[\ea\ea^T] & {\bf 0}\\
   {\bf 0} & \mbox{E}[\eb\eb^T] 
\end{array}\right] ,                        
\end{equation}
where the RHS matrix represents the true but unknown 
joint error covariance, which has zero cross covariance, 
$\ea\eb^T={\bf 0}$, due to the assumption
of statistical independence.  The estimated joint covariance
is a conservative estimate of the true joint covariance 
because in practice ${\bf A} \succeq \mbox{E}[\ea\ea^T]$ and
${\bf B} \succeq \mbox{E}[\eb\eb^T]$.  The latter inequalities
hold by design in that intentional efforts are made to ensure
that estimate error covariances do not underestimate
the actual squared errors associated with sensor and
kinematic models. 

Given a consistent joint covariance for two given
$n$-dimensional estimates, the Kalman filter defines the 
optimal linear projection of the $2n$-dimensional joint
estimate back to the $n$-dimensional state space of interest.
The result of the Kalman projection is a mean and covariance
estimate $({\bf c},{\bf C})$ that represents the optimal 
fusion of the two
given mean and covariance estimates.  In fact, if there
is no additional information available (e.g., distribution
information), then the Kalman fusion estimate is optimal
according to virtually any error criteria~\cite{maybeck1}. 

In the case of statistically independent estimates 
$({\bf a}_1,{\bf A}_1), ({\bf a}_2,{\bf A}_2), ...,
({\bf a}_m,{\bf A}_m)$, the Kalman fusion equations have
a particularly simple form~\cite{maybeck1}:
\begin{eqnarray}
{\bf C} & = & \left( \A_1^{-1} + \A_2^{-1} + ... 
                               + \A_m^{-1} \right)^{-1}\nonumber\\
{\bf c} & = & {\bf C} \left( \A_1^{-1}{\bf a}_1 + 
       \A_2^{-1}{\bf a}_2 + ... + \A_m^{-1}{\bf a}_m \right).
\end{eqnarray}
If its underlying assumptions hold (i.e., consistency and
independence), then the above Kalman equations ensure that
the fused estimate $({\bf c},{\bf C})$ is consistent, and
${\bf C} \preceq {\bf A}_i ~ \forall i, 1\leq i\leq m$.  
However, any presumption of statistical independence in 
practical data fusion contexts should be carefully 
considered.  Specifically, virtually any sensor is subject
to time-correlated errors induced by the particular
conditions of its use (e.g., changes in temperature, 
platform vibrations, relative humidity), and errors
associated with the nonlinear transformation of its
measurements (e.g., from local spherical 
coordinates to a global coordinate frame) are deterministic 
and therefore non-independent.

If estimates $({\bf a}_1,{\bf A}_1), ({\bf a}_2,{\bf A}_2), ...,
({\bf a}_m,{\bf A}_m)$ are each consistent, but not completely
independent, then it is possible for the Kalman fused estimate 
$({\bf c},{\bf C})$ to be {\em inconsistent}.  In fact, if
${\bf A}_i = \ea_i\ea_i^T~ \forall i, 1\leq i\leq m$, then
{\em any} degree of correlation guarantees inconsistency, i.e.,
${\bf C} \not\succeq \mbox{E}[\ec\ec^T]$.  The key point is that a
Kalman fused estimate is not guaranteed to be consistent 
even if each of its given estimates are consistent. The reason 
why the Kalman filter fails is because although two given
estimates $({\bf a}, {\bf A})$ and $({\bf b}, {\bf B})$ may be 
individually consistent, the implicit joint covariance may not be 
if independence is assumed when the cross covariance 
between the estimates is
${\bf X}=\mbox{E}[\ea\eb^T] \neq {\bf 0}$.  Specifically:
\begin{equation}
\left[\begin{array}{cc}
   \mbox{E}[\ea\ea^T] & {\bf 0}\\
   {\bf 0} & \mbox{E}[\eb\eb^T] 
\end{array}\right]               \not\succeq                           
\left[\begin{array}{cc}
   \mbox{E}[\ea\ea^T] & {\bf X}\\
   {\bf X}^T & \mbox{E}[\eb\eb^T] 
\end{array}\right]   .                      
\end{equation}  
In other words, the Kalman filter only fails to produce
a consistent fused estimate when the implicit {\em joint} 
estimate is inconsistent.  Although the equations are
very simple and elegant for independent estimates, the
Kalman filter is also defined generally for any consistent
joint covariance with ${\bf X} \neq {\bf 0}$~\cite{maybeck1}.  
Therefore, the Kalman 
filter is guaranteed to produce consistent estimates as 
long as the given estimates are consistent and their 
cross covariance is known~\cite{jaz}.
Unfortunately, this poses significant challenges.  
The first challenge is that the cross covariance 
information must in principle be determined
{\em exactly}.  This can be seen by examining the
difference between two joint covariance matrices with 
different cross terms:
\begin{equation}
\left[\begin{array}{cc}
   {\bf A} & {\bf X}\\
   {\bf X}^T & {\bf B}
\end{array}\right] -       
\left[\begin{array}{cc}
   {\bf A} & {\bf Y}\\
   {\bf Y}^T & {\bf B}
\end{array}\right]      =
\left[\begin{array}{cc}
   {\bf 0} & {\bf X}-{\bf Y}\\
   ({\bf X}-{\bf Y})^T & {\bf 0}
\end{array}\right] .         
\end{equation} 
The difference matrix is not PSD for any
case in which ${\bf X} \neq {\bf Y}$.  The need for 
absolutely perfect cross covariance information
presents difficulties when estimates are the products
of nonlinear operations (e.g., coordinate transformations,
kinematic time projections, human-derived estimates)
because the error processes are not 
perfectly modeled.  For example, the same approximate
nonlinear transformation equations may be applied to
convert different radar observations of an object to a
common coordinate frame, so the errors committed are
clearly not independent, but it may not be possible
to determine exact cross covariances.

Kalman's original derivation of his eponymous filter was 
based on orthogonal projection theory, and
the fact that there exists a simple Bayesian interpretation
of the result when error distributions are Gaussian was
presented only as an interesting special case~\cite{Kalman}.  
However, many subsequent references to the Kalman filter
incorrectly suggest that the Kalman filter {\em requires} assumptions
of Gaussianity\footnote{In fact, one commonly-cited motivation
for investigating the application of neural networks and
fuzzy logic is the claim that the Kalman filter imposes
restrictive Gaussianity assumptions that cannot be 
satisfied in many applications. The fact is that the use of
covariance upper bounds was recognized as necessary when the
first Kalman filters were implemented in the late 1960s, and
such bounds are incompatible with PDF interpretations. It was shown 
by Jazwinski in his classic 1970 book that the standard practice
of using covariance upper estimates does not undermine the integrity
of the Kalman filter\cite{jaz}.}. It turns out that 
the assumption of estimate independence is actually the
only problematic assumption because it typically cannot
be guaranteed in practice. Relaxing the independence assumption 
leads to the more
general fusion equations of the Covariance Intersection (CI)
method which has proven invaluable for a wide variety of practical
applications\cite{article:impactci,article:julierslam07}.

\section{Covariance Intersection (CI)}

The general mean and covariance data fusion problem 
can be formulated in terms of the joint covariance structure
that implicitly exists between a given pair of estimates
$({\bf a},{\bf A})$ and $({\bf b},{\bf B})$:
\begin{equation}
   \left[\begin{array}{cc}
   \A & {\bf X}\\
   {\bf X}^T & \B
   \end{array}\right]
\end{equation}
where ${\bf X}$ represents the actual, but unknown, 
cross covariance between the two estimates.
If ${\bf X}$ were known, then it would be possible to
apply more general formulations of the Kalman filter
equations to produce an optimal fused estimate.  Unfortunately,
these generalizations only guarantee consistency if the
cross covariance is known {\em exactly}, i.e., it cannot
be conservatively approximated in any way analogous to
the way conservative covariance estimates are
used.  

Without knowledge of ${\bf X}$, the only way to ensure 
consistency in the application of the Kalman filter is
to identify a joint covariance that is guaranteed to be
consistent based on the information available.  
In the present context,
therefore, a joint covariance matrix ${\bf M}$ must be 
determined such that:
\begin{equation}
   {\bf M}   ~ \succeq ~
   \left[\begin{array}{cc} 
   \A & {\bf X}\\
   {\bf X}^T & \B
   \end{array}\right]                       
\label{eqn:gendiff}
\end{equation}
\noindent for every possible cross covariance ${\bf X}$.  
It can be inferred from the symmetry of the 
unknown cross covariance 
information (i.e., ${\bf M}$ must be consistent
for any instantiation ${\bf X} = {\cal X}$ and
for ${\bf X} = -{\cal X}$)
that the off-diagonal blocks of 
${\bf M}$ should be 
zero, and its diagonal blocks must be sufficiently
larger than ${\bf A}$ and ${\bf B}$ to account
for the effects of all possible degrees of correlation 
among the error components of the mean estimates
${\bf a}$ and ${\bf b}$.  

It has been shown (appendix~14 of~\cite{phdthesis:uhlmann95}) that a
consistent and tight joint covariance ${\bf M}$ 
can be generated by selecting a scalar value
$\omega$, $0\leq\omega\leq 1$ as:
\begin{equation}
   \left[\begin{array}{cc} 
   \frac{1}{\omega}\A & {\bf 0}\\
   {\bf 0} & \frac{1}{(1-\omega)}\B
   \end{array}\right]            \succeq
   \left[\begin{array}{cc} 
   \A & {\bf X}\\
   {\bf X}^T & \B
   \end{array}\right]
\label{eqn:matci}
\end{equation}
where $\omega$ is chosen to minimize the
size (e.g., determinant) of the covariance
produced by the Kalman filter update 
equations for the estimates 
$({\bf a},{\bf \frac{1}{\omega}A})$ and 
$({\bf b},\frac{1}{1-\omega}{\bf B})$. This
can be generalized for an arbitrary number
of estimates:
\begin{equation}
\left[\begin{array}{cccc}
   \frac{1}{\omega_1}{\bf A}_1 & {\bf 0}     & ... & {\bf 0} \\
      {\bf 0} & \frac{1}{\omega_2}{\bf A}_2  & ... & {\bf 0} \\
      \vdots  & \vdots & \ddots & \vdots \\
      {\bf 0} & {\bf 0} & ... & \frac{1}{\omega_m}{\bf A}_m  
   \end{array}\right]  \succeq
\left[\begin{array}{cccc}
      {\bf A}_1 & {\bf X}_{1,2} & ... & {\bf X}_{1,m} \\
      {\bf X}_{2,2} & {\bf A}_2  & ... & {\bf X}_{2,m} \\
      \vdots  & \vdots & \ddots & \vdots \\
      {\bf X}_{m,1} & {\bf X}_{m,2} & ... & {\bf A}_m  
   \end{array}\right]               
\end{equation} 
where $\sum_{i=1}^{m} \omega_i = 1$, and the parameters
are determined to minimize the covariance resulting from
the fusion of the $n$ estimates. The above inequality provides
a general (and optimal) mechanism for obtaining a consistent
joint covariance when given only the block diagonals of an
unknown joint covariance matrix. The use of appropriate
$\omega$-paramterized covariances in the Kalman update
equations yields the CI fusion equations. 

The general approach of determining consistent joint 
covariances can also be applied to solve related problems.
For example, given estimates $({\bf a},{\bf A})$ and 
$({\bf b},{\bf B})$ that are correlated to an unknown
extent, the covariance of ${\bf a}+{\bf b}$
can be computed as $\frac{1}{\omega}\A + \frac{1}{1-\omega}\B$.
This is referred to as Covariance Addition (CA). 

\subsection{Covariance Union (CU)}

Covariance Intersection addresses the general form of the 
data fusion problem for mean and covariance estimates, but  
in practice a different problem can arise before data fusion 
can even be performed.  
Specifically, what is to be done if two estimates 
$({\bf a}, {\bf A})$ and $({\bf b}, {\bf B})$, purportedly 
relating to the state of the same real-world object, are 
determined to be mutually inconsistent with each other, 
i.e., the 
differences between their means is much larger than what 
can be expected based on their respective error covariance 
estimates?  For example, if two mean position estimates 
differ by more than a kilometer, but their respective 
covariances suggest that each mean is accurate to within 
a meter, then clearly something is wrong. 

One mechanism for detecting statistically significant 
deviations between estimates is to compute Mahalanobis
distances~\cite{mahal}.  The Mahalanobis distance 
between estimates $({\bf a}, {\bf A})$ and $({\bf b}, {\bf B})$ 
is defined as:
\begin{equation}
({\bf a} - {\bf b})^T \left({\bf A} + {\bf B}\right)^{-1}
({\bf a} - {\bf b}), 
\end{equation}
which is essentially just the squared distance between the
means as normalized by the sum of their respective covariances.
Intuitively, if the covariances are large, then a large difference
between the mean vectors ${\bf a}$ and ${\bf b}$ is not
surprising, so the Mahalanobis distance is small.  However,
if the covariances are very small, then even small differences
between the means may yield a large Mahalanobis distance.
A large Mahalanobis distance may tend to indicate that the 
estimates are not consistent with each other, but a user-defined 
threshold is required to define what constitutes an acceptable
deviation\footnote{It must be emphasized that the use of a
threshold on Mahalanobis distance is not the only possible mechanism 
for identifying potentially spurious estimates, but {\em some} 
user-defined mechanism is required. Otherwise there is no way
to distinguish fault conditions from low probability events. 
In other words, models for fault conditions are inherently
application-specific.}. When the threshold is exceeded, the
estimates are regarded as being contradictory and some kind
of action must be taken. Resolving such 
inconsistencies among estimates is sometimes referred to as 
{\em deconfliction}~\cite{cnis}.

The Covariance Intersection method guarantees consistency as 
long as the estimates to be fused are each consistent.  In 
the deconfliction problem it is only known that one of the 
estimates, either $({\bf a}, {\bf A})$ or 
$({\bf b}, {\bf B})$, is a consistent estimate of the state 
of the object of interest.  Because it is not generally 
possible to know which estimate is spurious, the only way to 
rigorously combine the estimates is to form a unioned 
estimate, $({\bf u}, {\bf U})$, that is guaranteed to be 
consistent with respect to {\em both} of the two estimates. Such 
a unioned estimate can be constructed\cite{article:uhlmann03} by computing a mean 
vector ${\bf u}$ and covariance matrix ${\bf U}$ such that:
\begin{eqnarray}
{\bf U} & \succeq & {\bf A} +
                 ({\bf u}-{\bf a})({\bf u}-{\bf a})^T\nonumber\\
{\bf U} & \succeq & {\bf B} +
                 ({\bf u}-{\bf b})({\bf u}-{\bf b})^T
\end{eqnarray}
where some measure of the size of ${\bf U}$, e.g., 
determinant, is minimized.  
This Covariance Union (CU) of the two estimates 
can be subsequently fused with other consistent estimates using 
CI. The CU equations generalize directly for the case of $m>2$ two estimates:
\begin{eqnarray}
{\bf U} & \succeq & {\bf A}_1 +
                 ({\bf u}-{\bf a}_1)({\bf u}-{\bf a}_1)^T\nonumber\\
{\bf U} & \succeq & {\bf A}_2 +
                 ({\bf u}-{\bf a}_2)({\bf u}-{\bf a}_2)^T\nonumber\\
~       & \vdots & ~\nonumber\\
{\bf U} & \succeq & {\bf A}_m +
                 ({\bf u}-{\bf a}_m)({\bf u}-{\bf a}_m)^T
\end{eqnarray}

Intuitively, the CU equations simply say that if the 
estimate $({\bf a}, {\bf A})$ is consistent, then the 
translation of the vector ${\bf a}$ to ${\bf u}$ will 
require its covariance to be enlarged by the addition of 
a matrix at least as large as the outer 
product of $({\bf u}-{\bf a})$ in order to be consistent.  
The same reasoning applies 
if the estimate $({\bf b}, {\bf B})$ is consistent.  
Covariance Union therefore determines the mean vector ${\bf u}$ 
having the smallest 
covariance ${\bf U}$ that is large enough to guarantee 
consistency regardless of which of the two given estimates 
is consistent.  The resulting covariance may be 
significantly larger than either of the given covariances, but 
this is an accurate 
reflection of the actual uncertainty that exists due to the 
conflict between the two estimates. The key fact is that
the CU estimate satisfies the definition of consistency.  

As a simple example of a CU construction, consider two estimates 
$({\bf a}, {\bf A})$ and $({\bf b}, {\bf B})$ of the
location of an object observed from two nodes in a network.
The estimate from the first node places the mean
position at ${\bf a} = [0, 0]^T$, and the second node places
it at ${\bf b} = [4, 4]^T$, and each has an error covariance equal 
to the identity matrix $\mbox{\bf I}$. If it is determined that the
two estimates are statistically inconsistent with each other, thus
implying that one of the estimates is not a consistent estimate
of the object's location, then deconfliction must be performed.
The optimal CU deconflicted estimate is: 
\begin{equation}
{\bf u} = \left[\begin{array}{c} 2 \\ 2 \end{array}\right]
 , \hspace*{0.5in} {\bf U} =
\left[\begin{array}{cc} 5 & 4\\ 4 & 5 \end{array} \right]
\end{equation}
It is straightforward to verify that this estimate $({\bf u}, 
{\bf U})$ is in fact consistent with respect to either/both of 
the estimates $({\bf a}, {\bf A})$ and $({\bf b}, {\bf B})$.
If $({\bf a}, {\bf A})$ is a consistent estimate of the target's 
state, then the covariance ${\bf U}$ for mean ${\bf u}$ must be 
greater than or equal to 
${\bf A} + ({\bf u}-{\bf a})({\bf u}-{\bf a})^T$,
which it is.
It can be verified that the estimate $({\bf u}, {\bf
U})$ is similarly consistent with respect to the estimate 
$({\bf b}, {\bf B})$.  Therefore, if either of the two estimates
represents a consistent estimate of the state of the object, 
then the CU estimate is also consistent.

\section{General Covariance Union (GCU)}

Consistency of the CU equations rests on an implicit assumption 
that the estimates to be combined are not correlated.  Specifically, the CU 
inequalities for the estimate $({\bf u},{\bf U})$:
\begin{eqnarray}
{\bf U} & \succeq  &{\bf A}_1 + ({\bf u}-{\bf a}_1)({\bf u}-{\bf a}_1)^T\nonumber\\
{\bf U} & \succeq & {\bf A}_2 + ({\bf u}-{\bf a}_2)({\bf u}-{\bf a}_2)^T\nonumber\\
     ~ & \vdots  &~\nonumber\\
{\bf U} & \succeq & {\bf A}_m + ({\bf u}-{\bf a}_m)({\bf u}-{\bf a}_m)^T\label{eqn:cu}
\end{eqnarray}
technically hold only if the errors associated with ${\bf A}_i$ 
are uncorrelated with ${\bf a}_i$. This is true for most types of common
process models, but not for certain kinds of recursive control and colored
noise models.

For example, suppose that we are given a $1D$ estimate
$(0,0)$, i.e., zero mean and zero variance. The underpinning assumptions of CU imply
that if the mean is translated to $1$, the new consistent estimate must have variance
$0+1^2=1$. If instead the mean is translated to $2$, the variance must be $0+2^2=4$. 
However, suppose
a 2-step translation is performed: First the estimate is translated $1$ unit to 
produce $(0+1,0+1^2)=(1,1)$, then that estimate is translated $1$ more unit. According to CU the
result should be $(1+1, 1+1^2)=(2,2)$, but this is clearly incorrect because translating
the estimate $(0,0)$ by $2$ units must produce $(2,0+2^2)=(2,4)$ to ensure consistency. 
The problem, of course, is that the translations in the sequence of steps are 
not independent.

More generally, given a mean and covariance estimate $({\bf a},{\bf A})$, the
covariance of ${\bf a}+{\bf x}$ is equal to 
${\bf A}+ {\bf x}{\bf x}^T$ if and only if the estimate error, $\tilde{{\bf a}}$, is 
independent of the translation ${\bf x}$. If they are correlated to an unknown extent,
then the covariance of ${\bf a}+{\bf x}$ must be formulated using Covariance Addition (CA)
as $ \frac{1}{\omega} {\bf A}+ \frac{1}{1-\omega}{\bf x}{\bf x}^T$. Applying CA
to the CU inequalities (\ref{eqn:cu}) gives the generalized GCU formulation:
\begin{eqnarray}
{\bf U} & \succeq  & \frac{1}{\omega_1} {\bf A}_1 + \frac{1}{1-\omega_1} ({\bf u}-{\bf a}_1)({\bf u}-{\bf a}_1)^T\nonumber\\
{\bf U} & \succeq  & \frac{1}{\omega_2} {\bf A}_2 + \frac{1}{1-\omega_2}({\bf u}-{\bf a}_2)({\bf u}-{\bf a}_2)^T\nonumber\\
     ~ & \vdots  &~\nonumber\\
{\bf U} & \succeq & \frac{1}{\omega_m} {\bf A}_m + \frac{1}{1-\omega_m}({\bf u}-{\bf a}_m)({\bf u}-{\bf a}_m)^T\label{eqn:gcu}
\end{eqnarray}
where the optimization problem now involves $n$ free parameters ($0\leq\omega_i\leq 1$), like CI, but with the semidefinite inequalities of standard CU. This appears to
represent a daunting challenge to solve efficiently. However, plots of solutions
for low-dimensional problems provide a critical insight. 
Specifically, Figure~1 depicts the 1-sigma solution contour that
minimally encloses the 1-sigma contours of the estimates that are unioned. This
suggests -- but of course does not prove -- that the GCU solution is equivalent
to the corresponding MEE solution obtained by interpreting the covariances as
ellipsoids defined by their 1-sigma contours. That they are in fact equivalent is 
established in the following section.

\begin{figure}
  \begin{center}
    \input{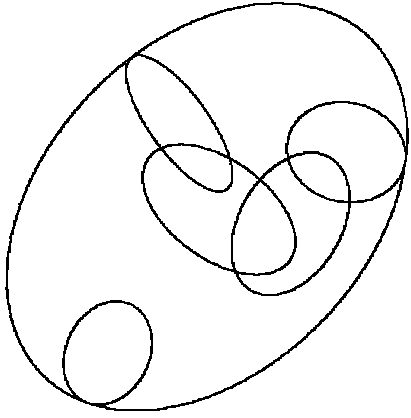}
  \end{center}
  \caption{\footnotesize A typical General Covariance Union (GCU) optimization problem. 
           The five small covariance ellipses represent the five estimates 
           to be unioned. Surprisingly, they seem to be circumscribed by the 
           minimal-determinant GCU solution.}
\end{figure}

\section{The Equivalence of GCU and MEE}
%**********************************************

This section provides a formal proof of the equivalence of MEE and 
GCU\footnote{To make the exposition concise our proof assumes that the covariance
matrices are nonsingular, but the equivalence holds generally.}
Our treatment of ellipsoid enclosure follows the discussion in \S3.7.1 
of~\cite{lmibook}, 
so several of its equations are repeated here for easy reference.

In~\cite{lmibook}, the enclosing and enclosed ellipsoids are represented 
via parameter triples $\ellipseTriple$ 
where index $0$ refers to the enclosing ellipsoid:
\begin{eqnarray}\label{eqn:BoydEllipsoids}
  \{ x \in \setRn \mid  x^T A_i x + 2 x^T b_i + c_i 
      \le 0 \}, \quad i = 0 \ldots m, \\
                               \nonumber\\
  A_i = A_i^T \succ 0 \quad \mbox{and} \quad 
      b_i^T A_i^{-1} b_i - c_i > 0\nonumber
\end{eqnarray}
The triples $\ellipseTriple$ are homogenous so the authors normalize 
the enclosing ellipsoid's parameters via 
$c_0 = b_0^T A_0^{-1} b_0 - 1$ 
and then use the S-lemma to represent the ellipsoid-enclosure 
constraint as a matrix inequality:
\begin{equation}\label{eqn:BoydCMI}
\left[\begin{array}{cc} 
  A_0    &  b_0 \\
  b_0^T  &  b_0^T A_0^{-1} b_0 - 1 \end{array} \right] - \tau_i 
\left[\begin{array}{cc}
  A_i    &   b_i\\
  b_i^T  &   c_i \end{array} \right] \, \preceq \, 0,
\quad \tau_i \ge 0, \quad i = 1 \ldots m.
\end{equation}
That inequality is nonlinear so they use a Schur complement argument to 
expand it:
\begin{equation}\label{eqn:BoydLMI}
\left[\begin{array}{ccc} 
  A_0    &  b_0   &  0\\
  b_0^T  &  - 1   &  b_0^T \\
  0      &  b_0   &  -A_0 \end{array} \right] - \tau_i 
\left[\begin{array}{ccc}
  A_i    &   b_i  &  0 \\
  b_i^T  &   c_i  &   0^{T}\\
  0      &   0    &  0 \end{array} \right] \, \preceq \, 0,
\quad \tau_i \ge 0, \quad i = 1 \ldots m.
\end{equation}
The minimum-determinant ellipsoid enclosure can then be posed as a 
Maxdet~\cite{maxdet} optimization:
\begin{equation}\label{eqn:BoydEllFormulation}
\begin{array}{cl}
  \mbox{minimize}   & \log \det A_0^{-1} \\
  \mbox{subject to} & A_0 \succ 0, \quad \tau_i \ge 0, \quad i = 1 \ldots m. \\
                    & \left[\begin{array}{ccc} 
  A_0    &  b_0   &  0\\
  b_0^T  &  - 1   &  b_0^T \\
  0      &  b_0   &  -A_0 \end{array} \right] - \tau_i 
\left[\begin{array}{ccc}
  A_i    &   b_i  &  0 \\
  b_i^T  &   c_i  &  0 \\
  0      &   0    &  0 \end{array} \right] \, \preceq \, 0, \quad i = 1 \ldots m.
\end{array}
\end{equation}

\subsection{Covariance ellipsoid enclosure as a Maxdet}
\label{sct:GCUEllipsoidFormulation}
%******************************************************

In this section, 
we translate the ellipse-enclosure formulation (\ref{eqn:BoydEllFormulation}) 
from the original $\{A,b,c\}$ triples to 
covariance/mean estimation parameters $(u,U)$.

The covariance ellipsoid $\covellipse{u}{U}$ for an estimate at 
mean $u$ with covariance matrix $U$ is defined as:
\begin{align}\label{eqn:CUEllipsoids}
\covellipse{u}{U} ~ = ~ & 
   \{ \, x \in \setRn \mid (x - u)^T U^{-1} (x - u) \le 1 \} \nonumber \\
= ~ & \{\, x \in \setRn \mid 
            x^T U^{-1} x + 2 x^T (-U^{-1} u) + (u^T U^{-1} u - 1) \le 0\}
\end{align}
Comparison of (\ref{eqn:BoydEllipsoids}) and (\ref{eqn:CUEllipsoids}) 
reveals the following correspondences: 
\begin{equation}\label{eqn:EllipseCorrespondence}
A \Leftrightarrow U^{-1}, \quad
b \Leftrightarrow - U^{-1} u, \quad
c \Leftrightarrow u^T U^{-1} u - 1
\end{equation}
Apply (\ref{eqn:EllipseCorrespondence}) to the nonlinear matrix inequality 
(\ref{eqn:BoydCMI}) - 
we will use the result in \S\ref{sct:singleEllipsoid}. 
The result is a GCU equation so $A_i$ now denotes a 
constraint covariance, and $a_i$ its mean: 
\begin{equation}\label{eqn:CuCMI}
\left[\begin{array}{cc} 
  U^{-1}          &  - U^{-1} u \\
  (- U^{-1} u)^T  &  u^T U^{-1} u - 1 \end{array} \right] - \tau_i 
\left[\begin{array}{cc}
  A_i^{-1}            &   - A_i^{-1} a_i \\
  (- A_i^{-1} a_i)^T  &   a_i^T A_i^{-1} a_i - 1 \end{array} \right] \, \preceq \, 0,
\quad \tau_i \ge 0, \quad i = 1 \ldots m.
\end{equation}
Returning to the covariance ellipsoid Maxdet, 
apply (\ref{eqn:EllipseCorrespondence}) to (\ref{eqn:BoydEllFormulation}): 
\begin{equation}\label{eqn:CuEllFormulation}
\begin{array}{cl}
  \mbox{minimize}    &  \log \det U \\
  \mbox{subject to}  &  U \succ 0, \quad \tau_i \ge 0, \quad i = 1 \ldots m,
\end{array}
\end{equation}
\begin{equation*}
 \left[ \begin{array}{ccc} 
  U^{-1}          &  - U^{-1} u  &  0 \\
  (- U^{-1} u)^T  &  - 1         &  (- U^{-1} u)^T \\ 
  0               &  - U^{-1} u  &  - U^{-1}  \end{array} \right] - \tau_i 
\left[\begin{array}{ccc}
  A_i^{-1}            &   - A_i^{-1} a_i          &  0 \\
  (- A_i^{-1} a_i)^T  &   a_i^T A_i^{-1} a_i - 1  &  0 \\
  0                   &   0                       &  0
\end{array} \right] \, \preceq \, 0,
\quad i = 1 \ldots m.
\end{equation*}
The matrix inequalities in (\ref{eqn:CuEllFormulation}) are nonlinear, 
due to the (quadratic) $U^{-1}u$ terms. 
Linearize them with change of variables $v = - U^{-1}u$. 
Also, let $W = U^{-1}$ to obtain a more standardized 
Maxdet formulation:
\begin{equation}\label{eqn:CuEllMaxdet}
\begin{array}{cl}
  \mbox{minimize}    &  \log \det W^{-1} \\
  \mbox{subject to}  &  W \succ 0, \quad \tau_i \ge 0, \quad i = 1 \ldots m, \\
                     &  \left[ \begin{array}{ccc} 
                           W      &  v      &  0 \\
                           v^T    &  - 1    &  v^T \\ 
                           0      &  v      &  - W \end{array} \right] - \tau_i 
\left[\begin{array}{ccc}
  A_i^{-1}            &   - A_i^{-1} a_i          &  0 \\
  (- A_i^{-1} a_i)^T  &   a_i^T A_i^{-1} a_i - 1  &  0 \\
  0                   &   0                       &  0
\end{array} \right] \preceq 0,
\quad i = 1 \ldots m.
\end{array}
\end{equation}
After solving (\ref{eqn:CuEllMaxdet}) for $v$ and $W$, 
the optimal $(u,U)$ are recovered as $U = W^{-1}$ and $u = -U v$.

\subsection{Enclosure of a single ellipsoid}
\label{sct:singleEllipsoid}
%*****************************************************************

Consider a GCU constraint from (\ref{eqn:gcu}) 
for a single estimate $(a,A)$ where $A$ is full rank:
\begin{equation}\label{eqn:GCUenclosure}
  U \, \succeq \, \frac{1}{\omega} A + 
           \frac{1}{1-\omega} (u - a)(u - a)^T, 
    \quad \omega \in [0,1]
\end{equation}
We shall prove that (\ref{eqn:GCUenclosure}) defines 
the set of all estimates $(u,U)$ whose covariance ellipsoid 
$\covellipse{u}{U}$ encloses $\covellipse{a}{A}$. 
It then follows that the combined GCU constraints 
(\ref{eqn:gcu}) 
define the set of all covariance ellipsoids which enclose all the 
$\covellipse{a_k}{A_k}$ 
so the solution to the (\ref{eqn:CuEllMaxdet}) Maxdet formulation is a 
minimum-determinant solution to (\ref{eqn:gcu}).

Note that (\ref{eqn:GCUenclosure}) only depends on the difference between 
$u$ and $a$ rather than their absolute values, 
and must hold for any arbitrary value of $u$. 
Therefore we may, without loss of generality, assume that $a = 0$ 
(simple coordinate shift) which simplifies 
(\ref{eqn:GCUenclosure}) to:
\begin{equation}\label{eqn:GCUorigin}
U \, \succeq \, \frac{1}{\omega} A + 
           \frac{1}{1-\omega} \, u u^T, 
    \quad \omega \in [0,1]
\end{equation}
A similar argument can be applied to a single instance from 
the ellipsoid-enclosure inequalities (\ref{eqn:CuCMI}): 
it must hold for arbitrary values of $u$, 
and ellipsoid enclosure is a geometric property unaffected by 
coordinate shifts. 
So we may again assume that $a = 0$, 
which simplifies (\ref{eqn:CuCMI}) to:
\begin{equation}\label{eqn:CuSingleCMI1}
\left[\begin{array}{cc} 
  U^{-1}            &  - U^{-1} u \\
  (- U^{-1} u)^T    &  u^T U^{-1} u - 1 \end{array} \right] - \tau 
\left[\begin{array}{cc}
  A^{-1}   &   0 \\
  0^T      &   - 1 \end{array} \right] \, \preceq \, 0, \quad \tau \ge 0
\end{equation}
Equation (\ref{eqn:CuSingleCMI1}) implies the following scalar inequality for 
the lower right-hand main diagonal entry:
$u^T U^{-1} u -1 + \tau \le 0$. 
Since $u^T U^{-1} u \ge 0$ we must have $\tau \le 1$. 
Make that restriction explicit, to bring $\tau$ 
in line with the equation (\ref{eqn:GCUorigin}) $\omega$ parameter:
\begin{equation}\label{eqn:CuSingleCMI}
\left[\begin{array}{cc} 
  U^{-1}            &  - U^{-1} u \\
  (- U^{-1} u)^T    &  u^T U^{-1} u - 1 \end{array} \right] - \tau 
\left[\begin{array}{cc}
  A^{-1}   &   0 \\
  0^T      &   - 1 \end{array} \right] \, \preceq \, 0, \quad \tau \in [0,1]
\end{equation}
We will prove equivalence by demonstrating that (\ref{eqn:GCUorigin}) and 
(\ref{eqn:CuSingleCMI}) are mutual inverses. 
But first we must consider the singular points. 
For (\ref{eqn:GCUorigin}) they are at $\omega = 0$ and $\omega = 1$.

As $\omega$ approaches 0, 
$U$ becomes unbounded and therefore $\covellipse{u}{U}$ 
encloses any (bounded) covariance ellipsoid $\covellipse{0}{A}$. 
The behavior for $\omega = 1$ depends on the value of $u$: 
if $u = 0$ then the ellipsoids are concentric and 
(\ref{eqn:GCUorigin}) reduces to the expected enclosure requirement 
$U \succeq A$.
If $u \neq 0$ then the right-hand side approaches the sum of $A$ and an 
unbounded rank-1 adjustment in the direction of $u$. 
The result is an elliptic cylinder - with radial axis through 
$0$ and $u$ - which tightly encloses $\covellipse{0}{A}$. 
So both singular points lead to enclosure of $\covellipse{0}{A}$.

Now that the singular points have been considered, 
we can simply invert (\ref{eqn:GCUorigin}): 
move the rank-1 term to the left-hand side and apply the 
Sherman-Morrison formula. 
The result simplifies to:
\begin{equation}\label{eqn:GCUoriginSherMorr} 
  U^{-1} +
    \frac{ (U^{-1} u) (U^{-1} u)^T } 
         { 1 - \omega - u^T U^{-1} u } \\
         \, \preceq \, \omega A^{-1}
\end{equation}
The rank-1 term in (\ref{eqn:GCUoriginSherMorr}) is 
bounded and positive semidefinite, 
since it is related to the bounded and positive semidefinite rank-1 term 
$\tfrac{1}{1-\omega} \, u u^T$ from (\ref{eqn:GCUorigin}). 
Therefore, $1 - \omega - u^T U^{-1} u > 0$ 
so we can apply a Schur complement:
\begin{equation}
\left[\begin{array}{cc}
  \omega A^{-1} - U^{-1}  &   U^{-1} u \\
  (U^{-1}u)^T             &   1 - \omega - u^T U^{-1} u 
      \end{array} \right] \, \succeq \, 0
\end{equation}
Finally, rearrange terms:
\begin{equation}\label{eqn:OmegaSingleCMI}
\left[\begin{array}{cc} 
  U^{-1}            &   - U^{-1} u \\
  (- U^{-1} u)^T    &   u^T U^{-1} u - 1 \end{array} \right] - \omega 
\left[\begin{array}{cc}
  A^{-1}   &   0 \\
  0^T      &   - 1 \end{array} \right] \, \preceq \, 0,
\quad \omega \in (0,1)
\end{equation}
Equation (\ref{eqn:OmegaSingleCMI}) has the same form as 
(\ref{eqn:CuSingleCMI}) 
(except for the singular points $\omega = 0$ and $\omega = 1$ 
which were considered earlier). 
This completes the proof.

\section{Discussion}

In this paper we have established a connection between covariance-based
and bounded-region models of uncertainty. In particular, we have 
shown that GCU and MEE have equivalent solutions. Thus, techniques
for solving one problem can be directly applied to solve instances
of the other. This equivalence is also suggestive of a potentially more
general framework that subsumes the covariance and bounded-region 
approaches.  

{\em Acknowledgements}: This work was supported by grants from the
Office of Naval Research and the Naval Research Laboratory.

  \bibliographystyle{plain}
  \bibliography{nrl-df}

\end{document}